\documentclass{amsart}
\usepackage{graphicx}
\usepackage{epsfig}

\usepackage{amsfonts}
\usepackage{isolatin1}

\title[{Why Delannoy numbers?}]{Why Delannoy numbers?}

\author[C. Banderier]{Cyril Banderier}
\email{Cyril.Banderier@lipn.univ-paris13.fr,\newline \tt http://www-lipn.univ-paris13.fr/$\sim$banderier}
\address{}
\author[S. Schwer]{Sylviane Schwer}
\email{Sylviane.Schwer@lipn.univ-paris13.fr,\newline \tt
http://www-lipn.univ-paris13.fr/$\sim$schwer/}
\address{\rm LIPN- UMR 7030
Université Paris Nord.\newline
99, avenue J.-B. Clément.
93430 Villetaneuse (France)}

\date{November 5, 2004}

\newcommand{\N}{{\mathbb N}}
\newcommand{\Z}{{\mathbb Z}}

\usepackage{multicol}

\newcounter{exampleno}
\begin{document}

\begin{abstract}
This article is not a research paper, 
but a little note on the history of combinatorics: 
We present here a tentative short biography
of Henri Delannoy, and a survey of his most notable works.
This answers to the question raised in the title,
as these works are related to lattice paths enumeration,
to the so-called Delannoy numbers,
and were the first general way to solve Ballot-like problems.
\end{abstract}
\maketitle

{\small \vskip-3mm
This version corresponds to an update (May 2002)
of the abstract submitted (February 2002) by the first author
to the 5th lattice path combinatorics
and discrete distributions conference (Athens, June 5-7, 2002).
The article will appear
in the  Journal of Statistical Planning and Inferences and this ArXiV
version has some minor typo corrections.}

\section{Classical lattice paths}

Before to tackle the question of 
Delannoy numbers and Delannoy lattice paths,
note that the classical number sequences
or lattice paths have the name of a 
mathematician: 
the Italian Leonardo Fibonacci ($\sim$1170--$\sim$1250),
the French Blaise Pascal (1623--1662),
the Swiss Jacob Bernoulli (1654--1705),
the Scottish James Stirling (1692--1770),
the Swiss Leonhard Euler (1707--1783),
the Belgian Eugène  Catalan (1814--1894),
the German Ernst Schr\"oder (1841--1902),
the German Walther von Dyck (1856--1934),
the Polish Jan {\L}ukasiewicz (1878--1956),
the American Eric Temple Bell (1883--1960),
the American Theodore Motzkin (1908-1970),
the Indian Tadepalli Venkata Narayana (1930--1987),
\ldots
It is quite amusing that some of them  are nowadays more
famous in combinatorics for problems which can be explained in terms of
  lattice paths  than in their original field (algebra
or logic for Dyck, Schr\"oder, and {\L}ukasiewicz\footnote{See the 
MacTutor History of Mathematics 
{\tt http://www-groups.dcs.st-and.ac.uk/$\sim$history/}}).

Fibonacci numbers appear in his 1202
{\em Liber abaci} (also spelled {\em abbaci})~\cite{Sigler}.
``Catalan numbers'' can be found in various works,
including~\cite{Cat,Segner}. Catalan called these numbers
``Segner numbers''; and the actual terminology is due to Netto
  who wrote the first classical introduction to combinatorics~\cite{Netto}.
The name ``Schr\"oder numbers'' honors the seminal paper~\cite{Sc70}
and can be found in  Comtet's ``Analyse
combinatoire''~\cite{Comtet70} and also in one of his
  articles published in 1970.
The name ``Motzkin  numbers'' can be  found
in~\cite{DoSh77} and is related to Motzkin's article~\cite{Motzkin}.
The name ``Narayana numbers'' was given by Kreweras by reference to the
article~\cite{Narayana1955} (these numbers were also independently
studied by John P. Runyon, a colleague of Riordan. These are called Runyon
numbers in Riordan's book~\cite{Riordan68}, p.17).
The name ``Dyck paths'' comes from the more usual ``Dyck words/Dyck Language''
which have been widely used for more than fifty years. We strongly recommend the lecture
of R. Stanley, which gives some comments about the surprisingly old
origin of these names and problems (cf pp. 212--213 of~\cite{Stanley}).

\pagebreak
\section{Delannoy numbers}
Delannoy is another ``famous'' name which is associated to an 
integer sequence related to lattice paths enumeration.
Delannoy's numbers indeed
correspond to the sequence $(D_{n,k})_{n,k\in\N}$, the number
of walks from $(0,0)$ to $(n,k)$, with
jumps $(0,1)$, $(1,1)$, or $(1,0)$.

\[\begin{array}{cccccccccc}
1&19&181&1159&5641  & 22363 & 75517 &  224143& 598417 & {\bf 1462563}
\\
1&17&145&833 &3649  & 13073 & 40081 & 108545 & {\bf 265729}& 598417
\\
1&15&113&575 &2241 & 7183  &  19825 & {\bf 48639} & 108545& 224143
\\
1&13&85 &377 &1289 &3653 & {\bf 8989} & 19825 & 40081& 75517
\\
1&11&61 &231 &681 &{\bf 1683} &3653 & 7183 & 13073& 22363 
\\
1&9 &41 &129 &{\bf 321} &681 &1289 &2241  &3649 & 5641
\\
1&7 &25 &{\bf 63}  &129 &231 & 377 & 575 & 833 & 1159 
\\
1&5 &{\bf 13}&25  &41  &61&85 &113 &145&181 
\\
1&{\bf 3}&5  &7   &9   &11&13&15 &17&19
\\
{\bf 1}&1 &1  &1   &1   & 1& 1& 1&1&1
\end{array}\]

In this array,
the lower left entry is $D_{0,0}=1$
and the upper right entry is $D_{10,10}= 8097453$.
Entry with coordinates $(n,k)$ gives the number of Delannoy walks
from $(0,0)$ to $(n,k)$. The three steps $(0,1)$, $(1,1)$, and $(1,0)$
being respectively encoded by $x$, $y$ and $xy$,
the generating function of Delannoy walks is
$$F(x,y,t)=\sum_{n\geq 0} (x+y+xy)^n t^n = \frac{1}{1-t (x+y+xy)}\,,$$
where $t$ encodes the length (number of jumps) of the walk.

The {\em central} Delannoy numbers $D_{n,n}$~({\tt EIS
1850}\footnote{This number refers to the wonderful On-Line Encyclopedia of
Integer Sequences, see {\tt http://www.research.att.com/$\sim$njas/sequences/}})
are in bold in the above array.
They have appeared for several problems:
properties of lattice and posets,
number of domino tilings of
the Aztec diamond of order $n$ augmented by an
additional row of length $2n$ in the middle~\cite{SaZe94},
alignments between DNA sequences~\cite{torres}\dots
The generating function of the central Delannoy numbers is
\begin{eqnarray*}
D(z)&:=&\sum_{n\geq 0} D_{n,n}
z^n=[x^0]   \, \frac{1}{1-(zx+z/x+z)} = \frac{1}{\sqrt{1-6z+z^2}}\\
&=& 1+ 3\, z+ 13\, z^2+63 \,z^3+321\, z^4+1683\, z^5+8989\, z^6+
48639\, z^7+ O(z^8)\,. 
\end{eqnarray*}
The notation $[x^n] F(x)$ stands for the coefficient of $x^n$ in the
Taylor expansion of $F(x)$ at $x=0$.
The square-root expression is obtained by a resultant or a residue computation
(this is classical for the diagonal of rational generating functions).
This closed form for $D(z)$ gives, by singularity analysis (see the
nice book~\cite{FlSe05}):
\begin{eqnarray*}
D_{n,n}&=& \frac{(3+2\sqrt 2)^n}{\sqrt \pi \sqrt{3\sqrt 2-4}}
\left( \frac{n^{-1/2}}{2}- \frac{23\, n^{-3/2}}{32 (8+3\sqrt 2)}+
\frac{2401 \,n^{-5/2}}{2048 (113+72\sqrt 2)}+O(n^{-7/2})\right)\\
&\approx&  5.82842709^n\,  \big(.57268163\, n^{-1/2}- .06724283 \,
n^{-3/2}+ .00625063\, n^{-5/2}+\dots\big)  \,.
\end{eqnarray*}
One has also
$\displaystyle D_{n,k}=\sum_{i=0}^n \binom{n}{i}\binom{k}{i} 2^i$.
Quite often, people note that there is a link
between Legendre polynomials and Delannoy numbers~\cite{Go58,La52,MoZa63},
and indeed $D_{n,n}=P_n(3)$, but this is not a very relevant link
as there is no ``natural'' combinatorial correspondence between Legendre polynomials
and these lattice paths. 
Comtet~\cite{Comtet70} showed that the coefficients of any algebraic
generating function satisfy a linear recurrence (which allows to compute them in linear time).
For $d_n:=D_{n,n}$, it leads to  $(n+2) d_{n+2} - (6n+9) d_{n+1}
+(n+1) d_{n} =0$.

For this kind of lattice paths with jumps $-1,0,+1$, 
one has links with continuous fraction~\cite{Fl80}, 
with determinants~\cite{GeVi},  with context free
grammars~\cite{LaYe90}... 
Numerous generalizations have been investigated:
walks in the quarter plane~\cite{Fa},
 multi-dimensional lattice paths~\cite{AuSc03,Kr03,SuDu04,GoNa69,HaMo76}.

It is classical in probability theory (see~\cite{BePi94} for some 
discussions with a combinatorial flavor)
and more precisely in the theory of Brownian motion to consider the
following constraints for lattice paths:
\begin{figure}[ht]
\small
\begin{center}\renewcommand{\tabcolsep}{3pt}
\begin{tabular}{|c|c|c|}
\hline
walks  & ending anywhere & ending in 0\\
\hline
  \begin{tabular}{c}  unconstrained\\(on~$\Z$) \end{tabular}
& \begin{tabular}{c}
	\quad { }\\ { \includegraphics[width=3.5cm]{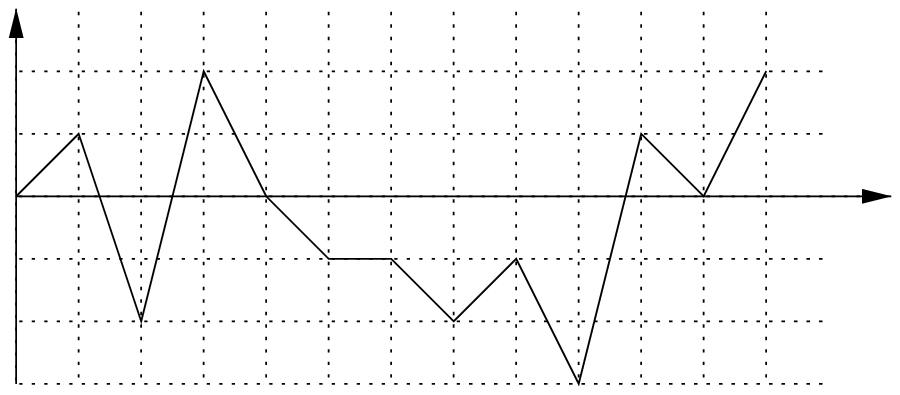}} \\
	walk ($\mathcal W$) \\
	\end{tabular}
& \begin{tabular}{c}
	\quad{ }\\ {\includegraphics[width=3.5cm]{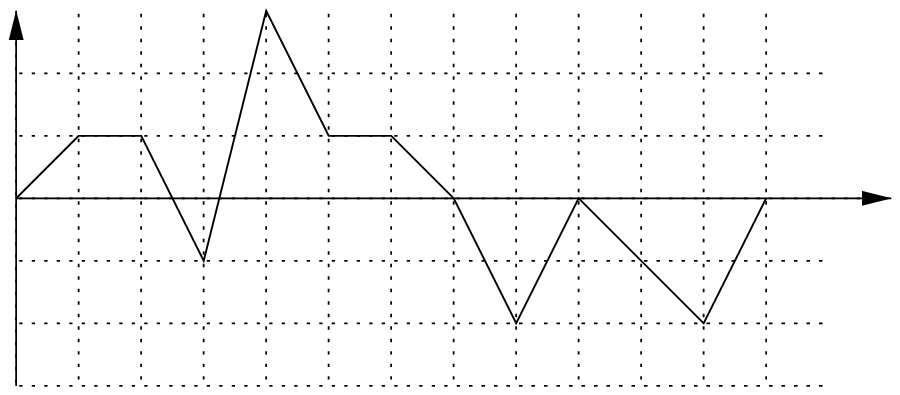}} \\
	bridge ($\mathcal B$)\\
	\end{tabular}	\\
\hline
\begin{tabular}{c} constrained\\ (on~$\N$) \end{tabular}
& \begin{tabular}{c}
	\quad { }\\ {\includegraphics[width=3.5cm]{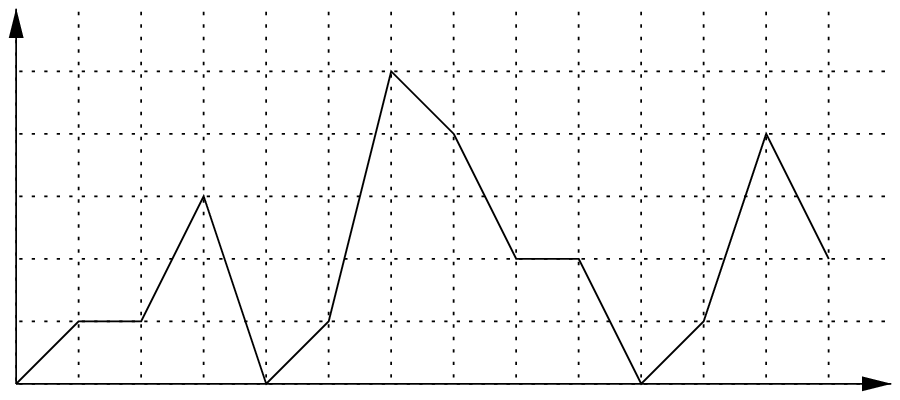}} \\
	meander ($\mathcal M$)\\
	\end{tabular}
&  \begin{tabular}{c}
	\quad { }\\	{\includegraphics[width=3.5cm]{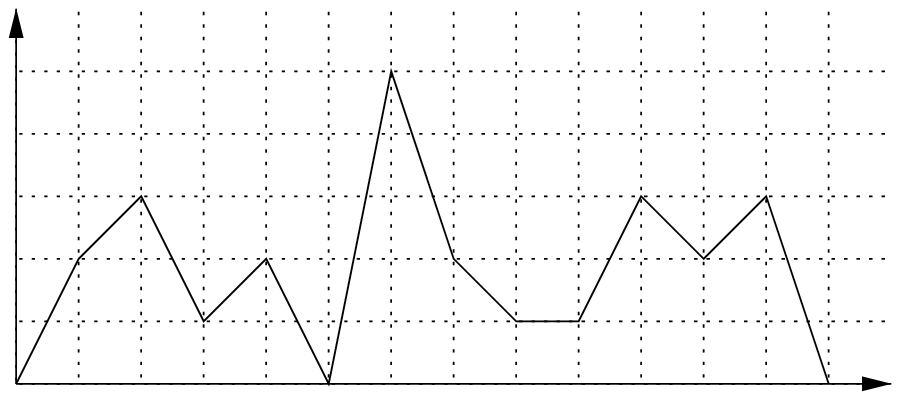}} \\
	excursion ($\mathcal E$)\\
	\end{tabular}\\
\hline
\end{tabular}
\end{center}

\vskip-2mm
\caption{\label{fig-4types}
	The four types of paths: walks, bridges, meanders, and
excursions. }
\end{figure}

\vskip-3mm
For these four kinds of walks and for any finite set of jumps,
there exists a nice formula for the corresponding generating function,
which appears to be algebraic and from which one can derive
the asymptotics and limit laws for several parameters of the lattice
paths~(see~\cite{BaFl01}).

Delannoy numbers $D_{n,n}$ correspond to bridges with a set of jumps
$\{+1,-1,0\}$ (where the 0 jump is in fact of length 2).
Consider now the language $\mathcal D$ of central  Delannoy paths, encoded
via the letters  $a, b, c$ (for the jumps $+1, -1, +0$ resp.).
Excursions with these jumps are called
Schr\"oder paths.
We note $\mathcal S$ the language of Schr\"oder paths (excursions) and 
$\bar {\mathcal S}$ 
the set of their mirror with respect to the $x$-axis.
Then, the natural combinatorial decomposition
$\mathcal D= (c^* a  \mathcal Sb+ c^* b\bar {\mathcal S} a)^* c^*$
(which means that one sees a Delannoy path [bridge] as a sequence of
Schr\"oder
paths [excursions] above or below the $x$-axis)  leads to
$$D(z^2)=\frac{1}{1-2z^2S(z) \frac{1}{1-z^2}}
\frac{1}{1-z^2}
\text{ where }
 S(z)=\frac{1-z^2-\sqrt{1-6z^2+z^4}}{2z^2}$$ is the
generating function of Schr\"oder paths.
This link between excursions and bridges is always easy to express
when the set of jumps is symmetric or with jumps of amplitude at most 1,
but there is also a relation between excursions and bridges
in a more general case (see~\cite{BaFl01} for  combinatorial and
analytical proofs).

Despite all these appearances of Delannoy numbers (see~\cite{Su03}
for a list of 29 objects counted by Delannoy numbers!),
the classical books in combinatorics or computer science
which are usually  accurate for ``redde Caesari quae sunt Caesaris''
  (e.g., Comtet, Stanley, Knuth)
are mute about this mysterious Delannoy.

\section{Henri Auguste Delannoy (1833-1915)}

\begin{figure}
\includegraphics[height=8cm,width=6cm]{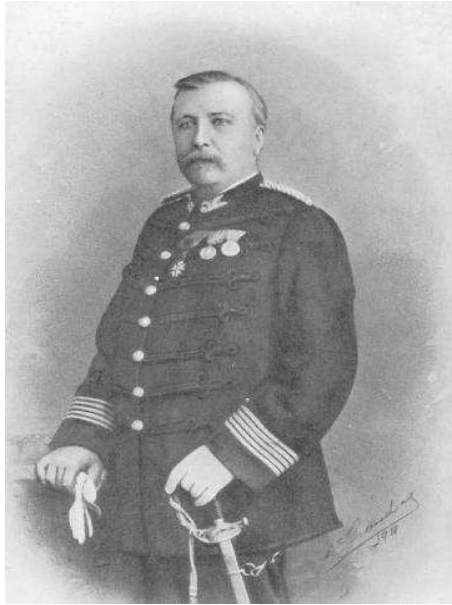}
\caption{Henri Auguste Delannoy (1833-1915). 
Portrait provided by
 the Société des Sciences Naturelles et Archéologiques de la Creuse,
where it is exhibited.}
\end{figure}

Some people suggested that ``Delannoy'' 
was related either to the French mathematician 
Charles Delaunay  (like in Delaunay triangulations)
or to the Russian mathematician
Boris Nikolaevich Delone, but this is not the case, as we shall see.

It is true that ``Delannoy'' sounds like (and actually is) 
a French family name
{\tt [d$\wedge$lanoa]} in approximative phonetic alphabet
  ({\tt $\wedge$} like in ``duck'' and {\tt a} in ``have'').
There are in fact thousands of Delannoy, mostly in the
North of France and in Belgium.
This toponym means ``de Lannoy'', that is to say
who originates from the town of Lannoy; 
``lannoy'' meaning a place with a lot of alders; yet how to find ``our'' Delannoy among all these homonyms?
The terminology ``Delannoy numbers'' became widely
used as it  can be found in Comtet's book
in the footnote from  exercise 20, p.93~\cite{Comtet70}:
``these numbers are often called Delannoy numbers''
without any reference.
In the English edition ``Advanced Combinatorics''~\cite{Comtet74},
the footnote becomes inserted in the text (p.81)
but there is still no reference.

In fact, it appears to be a good idea to look in Lucas' books 
(see~\cite{dec98} for some biographical informations on \'Edouard
Lucas [1842-1891]):
Indeed, in the second edition of the first volume of the {\em Récréations
mathématiques}~\cite{LucasRM91}, Lucas wrote in the
preface ``{\em J'adresse mes plus vifs remerciements \`a  mon ami sinc\`ere
et d\'evou\'e, Henry\footnote{There is no mistake, he was born Henry and
asked to change it into Henri. We use in this article the first name
Henri, as it was Delannoy's choice and as this was officially approved.}
  Delannoy,\ldots}'' and at page 13 of this
introduction we are told that Henri Delannoy was intendant. In the second
volume~\cite{LucasRM96}, the fourth recreation was dedicated to {\em
Monsieur Henri Delannoy, ancien
\'el\`eve de l'\'Ecole Polytechnique, sous-intendant militaire de Premi\`ere
classe}.  After the death of his friend Lucas in October 1891, Henri Delannoy contributed
with Lemoine and Laisant\footnote{The r\^ole played by each one is
explained in~\cite{AuDeSc03}.} 
to the publication of the third and fourth volumes\footnote{Available at the web site of the French National Library
{\tt http://gallica.bnf.fr/}}~\cite{LucasRM92,LucasRM94} 
as well as to the book {\em L'arithm\'etique
amusante}~\cite{Lucas95}.

Like most of the French military intendants, Delannoy was a student
from the \'Ecole Polytechnique
(which was the place where military officers received a scientific education).
From a database\footnote{Available at {\tt http://bibli.polytechnique.fr:4505/ALEPH0/}}
of the former students,  one knows that Henri (Auguste)
Delannoy is born in 1833 in Bourbonne-les-Bains (Haute-Marne, France).
His father was  Omère Benjamin Joseph Delannoy
(countable officer) and his mother was Françoise Delage;
they were living in the city of Bourges.
In 1853, he passed the \'Ecole Polytechnique
entrance exam  (with rank 62); then he graduated in 1854 with  rank 91/106
and finished with rank 67/94 in 1855.
It is quite funny that 
this database also contains, like for any other polytechnician,
a physical description of Henri Delannoy:
dark brown hair, average brow, average nose, blue eyes,
small mouth, round chin, round face, height: 1,68m!

In the archive center of the French Army (in the Château de Vincennes), one can
find his record under the number 61241. 
From this and~\cite{DicoU,Creuse}, we know that Delannoy was first in the
Artillery corps as sous-lieutenant (with rank 12/37 from the application school),
lieutenant (1857),
took part in the Italy campaign (27 May- 18 August 1859) and in the
Solférino battle (24 June 1859).
When he came back, 
he married his dulcinea Olympe-Marguerite Guillon
on the 10 November 1859. They had 2 daughters and one boy. 
Delannoy was promoted captain in 1863.  
He then became a supplier-administrator: Intendant-Adjoint in
1865, sous-intendant of third class in 1867, 
of second class in 1872, of first class in 1882 (he was now a widower). 
He spent three years  in Africa (6 Oct. 1866 - 25 Oct. 1869). 
He was the governor of
the military Hospital of Sidi-bel-Abes, Algeria, during the terrible typhus epidemic
(he belonged  the Supply Corps and they were in charge the
sanitary affairs).
He translated for himself and perhaps also for his hierarchy several
German books/notes about the 
Supply Corps. He took part in the 1870 war between France and Prussia.  
It is mentioned without explanation that he was in {\em Deutschland} 
on July 26, 1870 (that is, 4 days after the declaration of war\dots) 
and on  March 7, 1871  (that is, 3 days before the treaty of London\dots). 
He was decorated 
with the {\em m\'edaille d'Italie}, the {\em décoration sarde de la Valeur militaire},
 the {\em Croix de la Légion d'Honneur} on July 18, 1868, 
and the {\em Rosette d'Officier de la Légion d'Honneur} in December  20, 1886.
 He could have reached the highest military ranks, but he wanted peace
and decided to retire  (January 9, 1889)  
in Guéret (the main city of the French department ``la Creuse''),
beginning a second life dedicated to science and more particularly
to mathematics.

Looking for {\tt Delannoy*} in the Zentralblatt volumes of 1860--1920\footnote{
The ``Jahrbuch über die Fortschritte der Mathematik''
(annual review on the progresses of mathematics)
is available at {\tt http://www.emis.de/MATH/JFM/}},
gives the following nine
articles~\cite{Del88,Del89,Del90b,Del90,Del94,Del95,Del95b,Del97,Del98}.
The references, problems, methods, and solutions
used in these articles are similar to the ones  often mentioned in Lucas' books.
Most of Delannoy's articles
are signed by Monsieur (H.) Delannoy, 
military intendant in the city of Orléans (and later,
  retired military intendant in the city of Guéret). 
Delannoy was a quite active member of the French Mathematical Society (SMF)
in which he was admitted in 1882, introduced by Lucas and Laisant. He
  disappeared from the SMF's list in 1905, while he still contributed 
to l'Interm\'ediaire des Math\'ematiciens until 1910. 
This amateur mathematician then sank into oblivion
and we found no obituary in the SMF bulletins at the occasion of his
death on February 5, 1915.

In his death certificate, he was referred as president of the {\em Soci\'et\'e des
Sciences Naturelles et Archéologiques de la Creuse}.
This Society is in fact still very dynamic\footnote{{\tt
htpp://perso.wanadoo.fr/jp-l/SSC23/}}.
It eventually appears that this Society, over which Delannoy presided  
from 1896 to 1915, has some archives, a part of which was given
 by Delannoy's family. They include some biographies
written when Delannoy was still alive~\cite{GrAncy,DicoU}, a list of
his publications, and also 
an obituary and a short biography by members of the
Society~\cite{Creuse, Car76}.
We shall come back on Delannoy's works in this Society in Section~5 and
we now consider Delannoy's contribution to mathematics.

\section{Delannoy's mathematical work}
Delannoy began his mathematical life reading the mathematical recreations that
Lucas began to publish in 1879 in {\em La Revue Scientifique}. 
He was in contact with him in 1880 and began immediately to work with him, 
answering to letters of mathematicians transmitted by Lucas.

The first mention, in a mathematical work, to Delannoy
is in an article by Lucas ``Figurative arithmetics and permutations
(1883)''~\cite{Lucas83}, which deals with enumeration of
configurations of 8 queens-like problems (the simplest one being:
how to place $n$ tokens  on an $n\times n$ array, with no row or
column with 2 tokens). Delannoy is there credited for
having computed several sequences. 

Some years later, in 1886, Delannoy made his first 
mathematical public appearance in the annual meeting of  the
``Association Française pour l'Avancements des Sciences''. 
We now give the list of Delannoy's articles.

\subsection{Using a chessboard to solve arithmetical
problems(1886)~\cite{Del86}}

In this article, Delannoy comes back on Lucas' article  mentioned
above and explains how he can use a ``chessboard'' (in modern words:
an array) to get  the formula $\frac{n-k+1}{n+1} \binom{n}{n+k}$ for
the number of Dyck paths of length $n$
ending at altitude $k$ by using something which is not far
from what one calls now the  Desiré André
reflection principle, which was in fact published one year later~\cite{Andre87}.
Note that Feller says (without references, see pages 72 and 369
of~\cite{Fe71},  340 of~\cite{Fe68}) 
that Lord Kelvin ``method of images''
for solving some partial differential equations is a kind analytic equivalent
of the reflection principle in disguise.
However, if one has a look on William Thomson's letter to
Liouville~\cite{Ke45},  the link is rather mild and the
cleaver combinatorial ideas of André \& Delannoy cannot be attributed to Lord Kelvin.

Delannoy makes the link
$T_{x,y}=\binom{x}{x+y}-\binom{x-1}{x+y}=\frac{y-x+1}{y+1} \binom{x}{x+y}$
between entries from the rectangular array
(our walks on $\Z$, here given by the binomial coefficients)
and entries from the triangular array (walks constrained to remain in the
upper plane, our Dyck paths).
The numbers $T_{x,y}$ are called (in English)
``ballot numbers'', but they are also called Delannoy--Segner numbers
in Albert Sade's review (in the Mathematical Reviews)
of  Touchard's article~\cite{Touchard52}. 
Kreweras~\cite{Kreweras1992} and Penaud~\cite{Penaud95} 
follow this terminology (quoting Riordan or Errera~\cite{Errera}
but none of Delannoy's articles which all sank into oblivion).
In conclusion,  these ``Delannoy--(Segner)'' numbers $T_{x,y}$ 
are {\em not} the ``famous'' Delannoy numbers $D_{n,k}$ defined in Section 2.

\subsection{The length of the game (1888)~\cite{Del88}}

There are several  contributions of Rouché and Bertrand
in the Comptes Rendus de 
l'Académie des Sciences on the following problem  that they call {\em the
game}: ``two players have $n$ francs and play a game, at each round,
the winner gets one franc from his opponent. One stops when one of the
two players is ruined.'' When the game is fair, the probability to
be ruined at the beginning of the round $m$ is
(with $q=\frac{m-n}{2}$):
\begin{eqnarray*}
& & \frac{(-1)^{m-n}}{n} \sum_{k=1}^n (-1)^{k-1} \sin\big(\frac{(2k-1)
\pi}{2n}\big) \cos^{m-1}\big(\frac{(2k-1)\pi}{2n}\big)\\
& = & \frac{n}{2^{m-1}} \sum_{k=0}^{q/n} (-1)^k
\frac{2k+1}{\frac{m+n}{2}+kn} \binom{m-1}{q-kn}\,.
\end{eqnarray*}
Rouché proves the left hand part with some determinant computations
and Delannoy uses lattice paths to get the right hand part
(claiming justly that there was a mistake in Rouché's first formula).

One can see this problem as a Dyck walk in the strip $[-n,n]$,
that is why the formula is similar to the formula 14 from~\cite{Knuth72}
in their enumeration of planted plane trees of bounded height
(Feller~\cite{Fe68} gives also some comments on this).

\subsection{How to use a chessboard to solve various probability theory
problems (1889)~\cite{Del89}}
This is a potpourri of seven ballot-like or ruin-like problems
partially
solved by de Moivre, Laplace, Huyghens, Ampère, Rouché, Bertrand, André,~\ldots
for which  Delannoy presents his simple solutions,
obtained by his lattice paths enumeration method.
He calls the lattice ``chessboard''.
The different constraints corresponds to different kind of chessboards:
triangular for walks in the upper-plane, rectangular for unconstrained
walks, pentagonal for walks bounded from above,
hexagonal for walks in a strip (modern authors from statistical
physics sometimes talk about walks with a wall or two walls~\cite{Kr00}).
Delannoy numbers (and the two corresponding binomial formulae)
appear at page 51.
Delannoy says that it corresponds to the directed walk of a queen
(sic), and that this problem was  suggested to him by Laisant.
This (and the further advertisement by Lucas of Delannoy's works,
see e.g. p. 174 of~\cite{LucasTN91} on ``Delannoy's arithmetical
square'', which is exactly the array given in Section~2)
answers to the question raised in our title. The authors 
who later wrote about Delannoy numbers/arrays
then gave references to Lucas~\cite{LucasTN91}, whereas
 Delannoy's articles sank into oblivion.

\subsection{Various problems about the game (1890)~\cite{Del90}}
Using an enumeration argument, simplifying the sum that he obtained
  and then using the Stirling formula,
he gives the asymptotic result $\frac{1}{\sqrt {2\pi}} \sqrt{2n}$
as the difference between the number of won and lost games, after $2n$ games.
He also answers to other problems, {e.g.} what is the probability to
have a group of $2, 3, \ldots, 8$ cards of the same color in a packet of
32 cards.

\subsection{Formulae related the binomial coefficients (1890)~\cite{Del90b}}
He gives several binomial formulae,
such as $\sum_{k=0}^p (p-2k)^2\binom{p}{k} =p2^p$.

\subsection{On the geometrical trees and their use in the theory of
chemical compounds. (1894)~\cite{Del92,Del94}}

A chemist asked for some explanations of Cayley's results,
mentioned in a German review. Delannoy translated this
review and corrected a computational mistake,
giving his own method, without knowing~\cite{Ca57,Ca75,Ca81}.
This corresponds to the sequences {\tt EIS 22} (centered hydrocarbons
with $n$ atoms) and {\tt EIS 200}  (bicentered hydrocarbons with $n$
atoms).
Application of combinatorics to enumeration of chemical configurations
is a subject which will be later revisited by~Pólya~\cite{Po}.

\subsection{How to use a chessboard to solve some probability theory
problems (1895)~\cite{Del95}}
Delannoy makes a summary of 17 applications of his theory of
triangular/square/pentagonal/hexagonal chessboard.
The array of Delannoy numbers (see our Section~2) appears on the page 76
from this article.

\subsection{On a question of probabilities studied by d'Alembert
(1895)~\cite{Del95b}}
Delannoy corrects some mistakes
in Montfort's solution to a problem raised by d'Alembert.

\subsection{A question of undetermined analysis (1897)~\cite{Del97}}

A review (by Professor Lampe) of this article 
can be found in the Jahrbuch \"uber die Fortschritte der Mathematik.
However, we were not able to get this article. There were in fact two journals
whose name was ``Journal de Math\'ematiques élémentaires'' (one edited
by Vuibert and the other edited by Bourget/Longchamps),
neither of them seems to contain the quoted article.

\subsection{On the probability of simultaneous events (1898)~\cite{Del98}}
A priest
wrote an article
in which he was bravely contesting the "third Laplace principle"
$P(A \cap B)= P(A)P(B)$ for two independent events,
arguing with three examples. Delannoy shows that they present
a misunderstanding
of  "independent
events", which goes back to the original fuzzy definition by de Moivre.

\subsection{Contributions to ``L'Interm\'ediaire des
Math\'ematiciens''~\cite{Del08}}

This journal was created in 1894 by C.-A. Laisant and Émile Lemoine.
It is quite similar to the actual sci.math newsgroups.
This journal was indeed only made of problems/questions/solutions/answers.

During  the quoted period,
numerous famous mathematicians      made some contributions to this journal:
Appell, Borel, Brocard, Burali-Forti, Cantor, Catalan, Cayley, Cesàro,
Chebyshev,  Darboux, Dickson, Goursat, Hadamard, Hermite, Jumbert,
Hurwitz, Jensen, Jordan, Kempe,
Koenigs, Laisant, Landau, Laurent, Lemoine, Lerch, Lévy, Lindel{\"o}f,
Lipschitz, Moore, Nobel, Picard,  Rouché\dots

From 1894 until 1908 (date of his last mathematical contribution),
Delannoy was an active collaborator: he raised or solved around 70
questions/problems.
These are questions number 20, 29, 32, 51, 84, 95, 138, 139, 140,
141, 142, 155, 191, 192, 314, 330, 360, 371, 407, 424, 425, 443, 444,
451, 453, 493, 494, 514, 601, 602, 603, 664, 668, 749, 1090, 1304,
1360, 1459, 1471, 1479, 1551, 1552, 1578, 1659, 1723, 1869, 1875,
1894, 1922, 1925, 1926, 1938, 1939, 2074, 2076, 2077, 2091, 2195,
2212, 2216, 2251, 2305, 2325, 2452, 2455, 2583, 2638, 2648, 2868,
2873, 3326.

These contributions can be classified in three sets:
the problems and solutions related to
combinatorics (enumeration and applications to probabilistic
problems), problems and solutions related to
elementary number theory (representations of integers
  as sum of some powers, Fermat-like problems),
and questions/answers related to Lucas' books
(so mainly recreative mathematics,
but not so trivial problems as it includes, {\em e.g.}, the four
color problem).


\smallskip

To these articles,
perhaps one should add some {\em récréations} of~\cite{LucasRM94}
(compare the warning in its preface), and also some problems written by Lucas,
but with Delannoy's solutions. The same holds for articles written by
Lucas in {\em La Nature}.

\smallskip

Finally, there are some books~\cite{Catalan1892, Comtet74,Frolow1886,
LucasRM91, LucasTN91,   LucasRM94,Lucas95,LucasRM96} 
(Lucas, Frolow, and Catalan intensively corresponded with Delannoy for their books) 
or articles~\cite{Aeppli23,Bonin1993, Brenti95, ALS02,Jo92,
kaparthi1991,Kimperling01,Kreweras1992,Mo58,Paul82,
Peart00,Penaud95,Rouche1888,SaZe94,Schwer,Sulanke2000, Touchard52,
Traverso1917,Vassiliev94} 
  which mention either Delannoy numbers or
some of Delannoy's results/methods.


\section{Other Delannoy's works}

Besides mathematics, Delannoy painted watercolors  and, 
perhaps more importantly, studied history.
Indeed, from 1897 to 1914, he published 29 accurate
archaeological/historical articles in the  Mémoires de la Société
des Sciences Naturelles et archéologiques de la Creuse.

Let us give a taste of Delannoy's writer talent: here some titles
of his articles:
``On the signification of word {\em ieuru}'',
``One more word about {\em ieuru}'',
``A riot in Guéret in 1705'',
``Aubusson's tapestries'',
``A bigamist in Guéret'',
``Grapevines in the Creuse'',
a lot of studies ``Criminal trials in the Marche. The case \dots'',
several studies on abbeys
and some ``Critical list of the abbots from \dots'', and last but not least,
``An impotence trial in the 18th century''.
When he died, at the age of 81, about a dozen  other articles were
still in progress.

\bigskip
Delannoy is surely one of the last ``self-made'' mathematicians who
succeeded in getting a name in this field, rivaling professional
mathematicians.
What he discovered is nowadays well understood and can be classified
as ``basic enumerative combinatorics''. However, despite the
simplicity of his tools, it seems to us that Delannoy's work (and
more generally, the
underlying combinatorics) is a nice example of what could, but is
actually not taught to young students (or even in high-schools),
as an introduction to research in mathematics, also allowing the use of
computers and computer algebra softwares.
This kind of mathematics is only present at the mathematical
Olympiads. This attractive bridge between enumeration, geometry,
probability theory,
analysis, \dots deserves a better place.

It appears very clearly, thanks to the archives of the Society of Natural
Sciences and Archaeology, that besides his own publications, Delannoy played a great
r\^ole in checking proofs for numerous mathematicians and historians who wrote to have
his contribution~\cite{AuDeSc03}.
The archive from the Society and from Delannoy's family in Gu\'eret
reveals a true {\em honn\^ete homme}, as defined in the seventeenth century.

\newpage
\noindent{\bf Acknowledgements.}

The first author's interest to Delannoy numbers comes from a talk
that Marko Petkov{\v s}ek  gave
in the Algorithms Seminar at INRIA in 1999
(a summary of this talk can be found in~\cite{BaPe}).
As an example, he was dealing with chess king moves
(his general result about the nature of different multidimensional
recurrences can be found in the article~\cite{BMPe00}).
M. Petkov{\v s}ek asked the first author what he knew about 
Delannoy, and C. Banderier then started to conduct an investigation...

The second author's interest to Delannoy numbers comes from her own
works. As a researcher
in Temporal Representation and Reasoning, she developed a model
based on formal languages theory~\cite{Schwer97}
instead of the logical or relational approaches. This framework
allowed her to enumerate easily all
possible temporal relations between $n$ independent
events-chronologies. Then she tried to know if some of
these sequences were already known. In the $n=2$ case, Sloane's On-Line
EIS provided her the name of Delannoy. She asked everybody she met who
is Delannoy? She already started her own investigation 
when Philippe Flajolet sent her to the first author. 

\bigskip

Many people helped us in our dusty investigations,
we want to thank for their technical and friendly help:
Jean-Pierre Larduinat (webmaster), R\'egis Saint-James (secretary) and
\'Etienne Taillemitte (president)
from the Soci\'et\'e des Sciences Naturelles
et Arch\'eologiques de la Creuse (Guéret)\footnote{\tt http://perso.wanadoo.fr/jp-l/SSC23/},
Serge Paumier and family Desbaux, from Delannoy's family,
Muriel Colombier from the Archives D\'epartementales de la Creuse,
Guy Avizou, first maire-adjoint and vice-pr\'esident of the Conseil
R\'egional de la Creuse,
Silke Goebel from the Jahrbuch Project
(Electronic Research Archive for Mathematics, Karlsruhe)\footnote{\tt
http://www.emis.de/projects/JFM/}, Anja Becker from the library of the
Max-Planck-Institut für Informatik (Saarbrücken)\footnote{\tt
http://www.mpi-sb.mpg.de/services/library/}, Geneviève Deblock from
the Conservatoire
numérique des Arts et Métiers (Paris)\footnote{\tt http://cnum.cnam.fr/},
Brigitte Briot from the library of INRIA (Rocquencourt)\footnote{\tt
http://www-rocq.inria.fr/doc/},
  Céline Menil from the library of the
Maine university (Le Mans)\footnote{\tt http://bu.univ-lemans.fr/},
Solange Garnier, Francine Casas and Fran\c{c}oise Thierry from the 
library of the University of
Paris-Nord\footnote{http://www.univ-paris13.fr/},
  Nathalie Granottier from the Centre International de Rencontres Mathématiques
(Marseilles)\footnote{\tt
http://www.cirm.univ-mrs.fr/SitBib/Bibli/debut.html}, and the people
from
the library of Jussieu (Paris)\footnote{\tt http://bleuet.bius.jussieu.fr/}.

We also thank Jean-Michel Autebert, Philippe Flajolet, Nick Goldman,
and Peter John for their comments. We appreciated the careful reading by the referees.

\smallskip

During the realisation of this work, Cyril Banderier
 was partially supported by the Future and Emerging Technologies
   programme of the EU under contract number IST-1999-14186 (ALCOM-FT),
  by the INRIA postdoctoral programme and by the Max-Planck-Institut.

\newpage

Although very long, the following bibliography is not 
a complete list of the so vast literature on lattice paths, 
but it is mainly a tentative bibliography (by no way exhaustive)
about enumeration of lattice paths ``related'' to Delannoy lattice
paths and other Delannoy works.

\bibliographystyle{plain}
\bibliography{delannoy3}
 
\bigskip
\end{document}